\documentclass[journal,12pt,onecolumn,draftclsnofoot,]{IEEEtran}

\usepackage[cmex10]{amsmath}
\usepackage[utf8]{inputenc}
\usepackage{algorithm}
\usepackage{algpseudocode}
\usepackage{graphicx}

\usepackage[colorlinks]{hyperref}

\usepackage[cmex10]{amsmath}
\usepackage[utf8]{inputenc}
\usepackage{multirow}
\usepackage{soul}
\usepackage{algorithm}
\usepackage{algpseudocode}
\usepackage{graphicx}
\usepackage{dsfont}
\usepackage{xcolor}
\usepackage{cite}
\usepackage{tabularray}

\usepackage[utf8]{inputenc}
\usepackage{mathtools}
\usepackage[mathscr]{euscript}
\usepackage{amsmath}
\usepackage{amssymb}
\usepackage{amsfonts}
\usepackage{amssymb}
\usepackage{amsmath}
\usepackage{verbatim}
\usepackage[T1]{fontenc}
\usepackage[utf8]{inputenc}

\usepackage{latexsym}
\usepackage{enumerate}
\usepackage{indentfirst}
\usepackage{calligra}
\usepackage{ulem}
\usepackage{graphicx}
\usepackage{array}
\usepackage{float}
\usepackage{bbm}
\usepackage{mleftright}

\makeatletter
\let\old@ps@headings\ps@headings
\let\old@ps@IEEEtitlepagestyle\ps@IEEEtitlepagestyle
\def\psccfooter#1{%
    \def\ps@headings{%
        \old@ps@headings%
        \def\@oddfoot{\strut\hfill#1\hfill\strut}%
        \def\@evenfoot{\strut\hfill#1\hfill\strut}%
    }%
    \def\ps@IEEEtitlepagestyle{%
        \old@ps@IEEEtitlepagestyle%
        \def\@oddfoot{\strut\hfill#1\hfill\strut}%
        \def\@evenfoot{\strut\hfill#1\hfill\strut}%
    }%
    \ps@headings%
}
\makeatother
\usepackage{soul}

\usepackage{multirow}
\usepackage{soul,color}
\usepackage{graphicx}
\usepackage{dsfont}
\usepackage{subcaption}

\usepackage[utf8]{inputenc}
\usepackage{mathtools}
\usepackage[mathscr]{euscript}
\usepackage{amsmath}
\usepackage{amssymb}
\usepackage{amsfonts}
\usepackage{verbatim}
\usepackage[T1]{fontenc}
\usepackage[utf8]{inputenc}
\usepackage{latexsym}
\usepackage{enumerate}
\usepackage{indentfirst}
\usepackage{calligra}
\usepackage{ulem}

\usepackage{array}
\usepackage{float}
\usepackage{bbm}

\usepackage{eurosym}

\usepackage{hyperref}

\usepackage{caption}

\usepackage{tikz}
\usepackage{pgfplots}

\pgfplotsset{compat=1.8}
\usepgfplotslibrary{statistics}
\pgfmathdeclarefunction{fpumod}{2}{%
        \pgfmathfloatdivide{#1}{#2}%
        \pgfmathfloatint{\pgfmathresult}%
        \pgfmathfloatmultiply{\pgfmathresult}{#2}%
        \pgfmathfloatsubtract{#1}{\pgfmathresult}%
        \pgfmathfloatifapproxequalrel{\pgfmathresult}{#2}{\def\pgfmathresult{5}}{}%
    }

\usepackage{tabularx}
\usepackage{flushend}
\usepackage{textcomp}
\usepackage{xcolor}
\usepackage{import}

\usetikzlibrary{trees,decorations,shadows}
\tikzset{level 1/.style={sibling angle=45,level distance=4mm}}
\usetikzlibrary{arrows.meta}
\usepackage{forest}
\usetikzlibrary{external}

\let\oldtikzexternalgetnextfilename\tikzexternalgetnextfilename \renewcommand{\tikzexternalgetnextfilename}[1]{\oldtikzexternalgetnextfilename{#1}\expandafter\tikzsetnextfilename\expandafter{#1}}

\usepackage{pgfplotstable}
\usepackage[outline]{contour}
\contourlength{1.2pt}

\pgfplotsset{compat=1.13} 

\usetikzlibrary{spy}
\usetikzlibrary{calc}
\usetikzlibrary{fadings}
\usetikzlibrary{patterns}
\usetikzlibrary{shadows}
\usetikzlibrary{mindmap}
\usetikzlibrary{backgrounds}
\usetikzlibrary{shapes.symbols}
\usetikzlibrary{shapes.multipart}
\usetikzlibrary{shapes.geometric}
\usetikzlibrary{automata,positioning}
\usetikzlibrary{decorations.fractals} 
\usetikzlibrary{decorations.markings}
\usetikzlibrary{decorations.pathreplacing}
\usetikzlibrary{decorations.pathmorphing}
    
\usepackage{bm}

\usepackage{nomencl}
\makenomenclature
\tikzset{edge from parent/.style={segment angle=10,draw}}

\tikzset{
  my rounded corners/.append style={rounded corners=2pt},
}

\def\BibTeX{{\rm B\kern-.05em{\sc i\kern-.025em b}\kern-.08em
    T\kern-.1667em\lower.7ex\hbox{E}\kern-.125emX}}

\renewcommand{\nomgroup}[1]{%
     \ifthenelse{\equal{#1}{O}}{\item[\textit{Operators}]}{%
        \ifthenelse{\equal{#1}{I}}{\item[\textit{Indices}]}{%
            \ifthenelse{\equal{#1}{A}}{\item[\textit{Acronyms}]}{%
            `\ifthenelse{\equal{#1}{V}}{\item[\textit{Variables and parameters}]}{}}}}}
\usepackage{scalerel}
\usetikzlibrary{svg.path}
\definecolor{orcidlogocol}{HTML}{A6CE39}
\tikzset{
    orcidlogo/.pic={
        \fill[orcidlogocol] svg{M256,128c0,70.7-57.3,128-128,128C57.3,256,0,198.7,0,128C0,57.3,57.3,0,128,0C198.7,0,256,57.3,256,128z};
        \fill[white] svg{M86.3,186.2H70.9V79.1h15.4v48.4V186.2z}
        svg{M108.9,79.1h41.6c39.6,0,57,28.3,57,53.6c0,27.5-21.5,53.6-56.8,53.6h-41.8V79.1z M124.3,172.4h24.5c34.9,0,42.9-26.5,42.9-39.7c0-21.5-13.7-39.7-43.7-39.7h-23.7V172.4z}
        svg{M88.7,56.8c0,5.5-4.5,10.1-10.1,10.1c-5.6,0-10.1-4.6-10.1-10.1c0-5.6,4.5-10.1,10.1-10.1C84.2,46.7,88.7,51.3,88.7,56.8z};
    }
}

\newcommand\orcidicon[1]{\href{https://orcid.org/#1}{\mbox{\scalerel*{ \begin{tikzpicture}[yscale=-1,transform shape]
                \pic{orcidlogo};
                \end{tikzpicture}
            }{|}}}}

\allowdisplaybreaks

\usepackage{bm}
\usepackage{mdframed}
\usepackage{lipsum}

\newmdenv[leftline=false,rightline=false,linewidth=1pt]{topbot}
\begin{document}

\title{\huge{Multi-market Optimal Energy Storage Arbitrage with Capacity Blocking for Emergency Services}}

\author{
\IEEEauthorblockN{Md Umar Hashmi, Stephen Hardy, Dirk Van Hertem\\}
\IEEEauthorblockA{\small{\textit{KU Leuven}, division Electa \& \textit{EnergyVille}, 
Genk, Belgium,\\
(mdumar.hashmi, stephen.hardy, dirk.vanhertem)~@kuleuven.be\\ \vspace{19pt}}}
\and
\IEEEauthorblockN{\large{Harsha Nagarajan}\\}
\IEEEauthorblockA{\small{Applied Mathematics and Plasma Physics (T-5)\\
 Los Alamos National Laboratory, NM, USA, \\ 
 harsha@lanl.gov}}
}


 

\maketitle
\begin{abstract}
The future power system is increasingly interconnected via both AC and DC interconnectors. These interconnectors establish links between previously decoupled energy markets. In this paper, we propose an optimal multi-market energy storage arbitrage model that includes emergency service provisions for system operator(s). The model considers battery ramping and capacity constraints and utilizes operating envelopes calculated based on interconnector capacity, efficiency, dynamic energy injection and offshore wind generation in the day-ahead market.
The arbitrage model considers two separate electricity prices {for buying and selling of electricity in the two regions, connected via an interconnector}. 
Using disjunctive linearization of nonlinear terms, we exactly reformulate the inter-regional energy
arbitrage optimization as a mixed integer linear programming problem.
{We propose two capacity limit selection models for storage owners providing emergency services. The numerical analyses focus on two interconnections linking Belgium and the UK. The results are assessed based on revenue, operational cycles, payback period, shelf life and computation times.}

\end{abstract}

\pagebreak

\tableofcontents

\pagebreak

\section{Introduction}
With the large-scale adoption of AC and DC interconnections, the power network worldwide is becoming more and more interconnected, creating a new business case for grid-scale energy storage batteries. 
With the greater amount of interconnections, batteries can now simultaneously participate in more than one energy market and increase their revenue, compared to participating in just one system.
Prior use cases for grid-scale batteries \cite{link2vse, link3arena, bowen2019grid, schoenung2017green} indicate that they are often used for more than one application. One common responsibility of grid-scale batteries is providing emergency services under extreme operating conditions. Such obligations are often contractual and are defined in a {bilateral contract or power purchasing agreement (PPA). They may also be enforced via grid code.}

The goal of this work is to develop a computationally efficient optimization model for facilitating grid-scale battery participation in more than one energy market. {Further, we present models for determining the allocation of battery capacity to dedicate to emergency services. It's important to note that grid-scale batteries, often operate under bilateral agreements, which may require the provision of emergency services.} Thus, capacity blocking is crucial for assessing their true potential.


There exists limited literature on multi-regional energy market participation for storage.
Authors in \cite{hurta2022impact} observed that the bidding zone separation for the German and Austrian day-ahead power markets led to the decoupling of wholesale electricity prices.
They explored a scenario where a battery performs energy arbitrage in both German and Austrian day-ahead power markets operated by the Energy Exchange of Austria.
Prior works, \cite{bunn2010inefficient, soonee2006novel} discussed the benefits of inter-regional energy arbitrage or flows.
The authors in \cite{soonee2006novel} showcase that the energy trade between the northern and the southern grids of India, connected via HVDC back-to-back stations, led to {around 85 million euros} of savings for the years 2003 to 2005.
Authors in \cite{bunn2010inefficient} highlight the inefficiencies in the context of energy arbitrage in inter-regional electricity transmission, considering the example of the Anglo-French interconnector (IFA).

\subsection{Use cases of grid-scale batteries}
\label{sec:usecases}
{In this part, we elaborate on several use cases involving grid-scale batteries. These cases reveal that batteries often serve multiple objectives, including provisions for emergency services.}

    \textit{Victorian Big Battery}: A case resembling the system under discussion in this paper is exemplified in Australia.
    The 170 MW Victoria-New South Wales interconnector, is co-owned by AEMO and Transgrid. 
    On the Victorian side, there is a 300 MW, 450 MWh battery situated in Geelong, owned by Neoen. 
The battery utilization objectives are \cite{link2vse}:\\
{
$\bullet$ 
\textit{Energy Market Participation}: the battery engages in the energy market year-round, with 50 MW allocated for summer (Nov. to March) and the entire 300 MW available otherwise. \\$\bullet$ \textit{Capacity blocking for Emergency Services}: To mitigate the risk of unanticipated power cuts during the summer season, the battery blocks 250\,MW to act as a "virtual transmission line" between Victoria and New South Wales. This is sufficient to power over 1 million Victorian households for 30 minutes.}


\textit{The Hornsdale Power Reserve:}
at 100 MW/129 MWh, has operated in South Australia since December 2017. {The battery utilization objectives are:\\
{$\bullet$ \textit{Energy Market Participation}: 30 MW and 119 MWh capacity are bid directly into the market.\\$\bullet$ \textit{Emergency services}: 
It provides backup power during grid outages, ensuring stable frequency until slow generators start.
} This was first demonstrated in 2017 when the battery rapidly injected power into the grid, preventing a potential cascading blackout after a large coal plant tripped \cite{link3arena, bowen2019grid}.}

{
\textit{Moss Landing battery} was commissioned in Dec 2020 and is owned by Vistra Energy. The capacity of this grid-scale battery is 400MW and 1600 MWh. In addition to energy market participation, it is bounded by a 10-year resource adequacy contract for 100 MW/ 400 MWh with PG\&E \cite{mosslanding}.
}

{In contrast to the above examples, Manatee grid-scale battery (409 MW, 900 MWh) is owned and operated by Florida Power and Light (FPL). This utility-owned and operated battery aims to increase reliability during emergency situations like hurricanes.}


{As observed in the above use cases, the operational objectives of the battery may change based on the ownership. In this work, we assume the battery and the interconnector are owned by different entities. This is also in line with EU regulation 2019/943 \cite{euregulation}.}

\vspace{-5pt}
\subsection{Contributions}
The key contributions of the paper are:\\
$\bullet$ \textit{Optimization formulation}: A novel Mixed-Integer Linear Programming (MILP) based formulation for energy storage, performing energy arbitrage in more than one energy market simultaneously. The optimization formulation considers interconnector capacity limits, storage efficiency, interconnector flows and efficiency and also ensures that under no case the battery is charging and discharging simultaneously.\\
{$\bullet$ \textit{Energy Storage Profitability Assessment}: We provide a comprehensive framework for assessing the profitability of energy storage, including considerations for reserving a fraction of battery capacity for emergency services. This assessment accounts for the battery's operational and shelf life limitations while projecting annual simulated revenue.}\\
$\bullet$ \textit{Realistic case studies} are performed for (a) \href{https://www.nemolink.co.uk/}{NEMO Link} (case 1), and (b) the set of interconnections between Belgium and the UK connecting the Princess Elisabeth Energy Island (PEEI) in the North Sea (case 2). In these case studies, we (i) perform the sensitivity of interconnector rent on the performance indices, (ii) provide a solution to how much battery capacity can be blocked for emergency services, and (iii) answer whether the battery owners should consider reserving interconnector capacity for inter-regional energy trades or not. 

This paper is organized as follows.
Section \ref{section2}, describes the energy market models for test cases.
Section \ref{section3}, details the inter-regional energy arbitrage model as an MILP formulation.
Section \ref{section4} extends the inter-regional arbitrage model to consider the interconnector flows and capacity blocking.
Section \ref{section5} details the performance indices used for evaluating the simulations.
Section \ref{section6} describes the two case studies for interconnections connecting Belgium and the UK.
Section \ref{section7} concludes the paper.

\pagebreak

\section{Energy market modelling}
\label{section2}
\begin{figure}[!htbp]
	\center
	\includegraphics[width=6.5in]{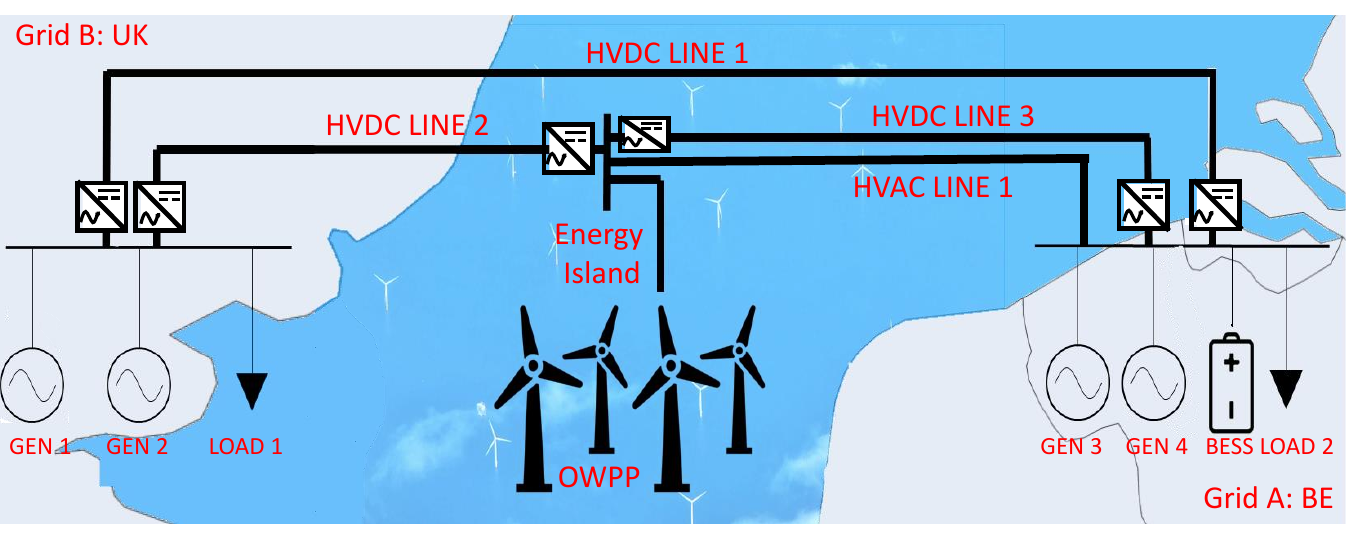}
	\vspace{-3pt}
	\caption{Stylized grids connected via an AC or DC interconnector with an offshore wind injection located in the centre.}
	\label{topo_fig}
\end{figure}

In this section, we outline a data-driven model for the Belgian-UK energy markets to assess the potential advantages of an inter-regional arbitrage strategy. The analysis begins by considering a single interconnector and then introduces a concept known as the Hybrid Offshore Asset (HOA), which combines offshore wind power and an interconnector as a shared resource. We explore two distinct case studies based on the Belgian offshore development plan up to 2030 as presented in the Federal Development Plan of the Belgian transmission system 2024-2034 \cite{elia_dev_plan}. The offshore grid considered is visualized in Fig. \ref{topo_fig}. The capacities of the transmission lines as well as which lines are considered in which test case are detailed in Table \ref{topo_table}.



\begin{table}[!htbp]
    \centering
        \caption{Test case transmission infrastructure.}
    \begin{tabular}{cccccc}
         Name&Transmission Line&from&to&GVA&test case\\\hline
         NEMO link&HVDC LINE 1&UK&BE&1.0&1\\ 
         Nautilus-UK&HVDC LINE 2&UK&PEEI&1.4&2\\
         Nautilus-BE&HVDC LINE 3&PEEI&BE&1.4&2\\
        -&HVAC LINE 1&PEEI&BE&2.1&2\\
        \end{tabular}
    \label{topo_table}
\end{table}
\subsection{Test Case 1: NEMO Link}
The first case focuses on the existing infrastructure, i.e., the 1 GW NEMO link interconnector connecting Belgium and the UK. The offshore wind generation is modelled as part of the home market, thus it is lumped in with other Belgian generation assets. Historical data from the years 2019 and 2020, sourced from the ENTSO-E transparency platform, are used to model the energy markets in Zone A (Belgium) and Zone B (UK). These years were selected due to the UK's departure from the common energy market at the beginning of 2020, known as "Brexit". The aim is to examine whether this event significantly impacts the feasibility of storage participating in inter-regional arbitrage.

The generators in zones A and B, denoted as $g_1$ and $g_3$, are assumed to be infinite sources with per-unit production costs represented by $P_{g_3,i}^{\text{\sc a}} \in \Psi^{\text{\sc p:a}}$ and $P_{g_1,i}^{\text{\sc b}} \in \Psi^{\text{\sc p:b}}$, where $\Psi^{\text{\sc p:a}}$ and $\Psi^{\text{\sc p:b}}$ are the historical day-ahead market clearing price time series for Belgium and the UK \cite{market_clearing}, respectively. 

Each market zone also includes a demand component with hourly values denoted as $D_{l_2,i}^{\text{\sc a}} \in \Psi^{\text{\sc d:a}}$ and $D_{l_1,i}^{\text{\sc b}} \in \Psi^{\text{\sc d:b}}$, where $\Psi^{\text{\sc d:a}}$ and $\Psi^{\text{\sc d:b}}$ correspond to the demand time series for Belgium and the UK \cite{trans_demand}, respectively. Furthermore, interconnector flow at each hour is defined as $L_{i}^{\text{\sc ab}} \in \Psi^{\text{\sc l:ab}}$, where $\Psi^{\text{\sc l:ab}}$ is the time series data of the net physical flows through the Nemo link from Belgium to the UK \cite{realized_flows}.



\subsection{Test Case 2: NEMO Link and the PEEI HOA}
In the second case, we introduce the 3.5 GW PEEI HOA, an energy island connected to both the UK and Belgium via the Nautilus link, a 1.4 GW HVDC interconnector. Additionally, the energy island is linked to Belgium through a 2.1 GW HVAC connection. 
To estimate energy production at the PEEI Offshore Wind Power Plant (OWPP), we employ the CorRES software \cite{Koivisto2019}, which utilizes historical meteorological data from the WRF database \cite{WRF_data} and introduces stochastic fluctuations and forecast errors over top to capture the inherent variability and uncertainties in renewable energy generation. {To ensure temporal compatibility with the existing data, simulations for years 2019 and 2020 are performed in CorRES.}

With the addition of the PEEI HOA, historical data alone becomes insufficient to model the energy markets. To address this, we introduce additional generators ($g_2$ and $g_4$) to each zone, as illustrated in Fig. \ref{topo_fig}. The hourly peak production of these assets is calculated based on the difference between hourly demand and a fixed "block size" parameter, assumed to be 1 GW. This parameter models the step-wise addition or removal of gas turbine-generating blocks from the energy market mix \cite{block_size}. As gas turbines often represent the marginal generating unit \cite{gas_is_marginal}, a reduction in generation of one block size due to offshore wind generation is assumed to shift the marginal pricing generator from gas to wind. 

The price of offshore wind results from an assumed bidding strategy of offering at 95\% of the Belgian energy market price. This idealized approach results in minimal curtailment for the OWPP while ensuring a healthy profit even when the OWPP becomes the marginal price-setting unit. The hourly pricing of generators $g_2$ and $g_4$ are therefore set at $P_{g_4,i}^{\text{\sc a}}=0.95 P_{g_3,i}^{\text{\sc a}} \in \Psi^{\text{\sc p:a}}$ and $P_{g_2,i}^{\text{\sc b}}=0.95 P_{g_1,i}^{\text{\sc b}} \in \Psi^{\text{\sc p:b}}$. 
Under these conditions, we calculate energy market prices by solving the optimal dispatch problem. This simulation provides estimates of both the market clearing prices and interconnector power flows after the addition of the PEEI. The market clearing prices for Belgium, the UK and the PEEI are determined from the dual variables ($\lambda_{i}^{\text{\sc be}}$, $\lambda_{i}^{\text{\sc ei}}$ and $\lambda_{i}^{\text{\sc uk}}$) of the nodal power balance equations, often referred to as the shadow prices.

\subsection{Flow-based energy market coupling}
Since the early 2010s, the EU has been working towards single-day-ahead market coupling across the EU via a gradual implementation of flow-based market coupling. Flow-based market coupling takes into account the physical constraints of power flows in interconnected grids, allowing for more efficient and reliable energy trading between different regions. It aims to maximize market efficiency and minimize congestion-related issues in the electricity grid \cite{entsoE_flow_based}. 

On January 31, 2020, the UK formally exited the European Union's common energy market, a development commonly referred to as "Brexit". Following Brexit, a provisional agreement governing energy markets between the EU and the UK went into effect. The provisional agreement is valid until 2026 \cite{UK_EU_market_agreement}. This agreement essentially maintained the pre-Brexit status quo, inhibiting any further progress towards the objective of flow-based market coupling.

While the energy market landscape beyond 2026 remains uncertain, there is optimism that an agreement can be reached, and the transition towards flow-based market coupling will resume. Consequently, in this study, we operate on the assumption that full implementation of flow-based market coupling will be achieved by the anticipated completion date of the PEEI project in 2030. This assumption carries the implication of improved and more efficient utilization of interconnector capacity in the second test case
{in the numerical results}
compared to the first.

\pagebreak

\section{Battery for inter-regional arbitrage}
\label{section3}
Given the grid configuration, a Battery Energy Storage System (BESS) can be strategically located at one end to capitalize on pricing disparities between the two grids during charge and discharge cycles. Furthermore, this battery can serve additional functions, including ancillary services and enhancing grid reliability. For the purposes of this study, we assume that the battery is located on the Belgian side of the interconnector and engages exclusively in inter-regional energy arbitrage.
{The battery is assumed to be a price-taker of electricity.}

\subsection{Battery model}
The battery model considers the ramping constraint, and the capacity constraint along with charging and discharging efficiencies denoted 
by $\eta_{\text{ch}}, \eta_{\text{dis}} \in (0,1]$, respectively. 
The energy optimization considers a change in energy levels of the battery {at time $i$} denoted as $x_i$. 
Change in battery energy level at $i$ is defined as $x_i = h \delta_i$, where $\delta_i \in [\delta_{\min}, \delta_{max}]$ $\forall i$ denotes ramp rate of the battery. 
$h$ denotes the sampling period.
$\delta_i> 0 $ when the battery is charging and vice versa. Note, $\delta_i$ is in units of power (MW) and $x_i$ is in units of energy (MWh). 
The battery charge level is denoted as
\begin{equation}
b_i = b_{i-1} + x_i, \quad b_i\in [b_{\min},b_{\max}], \forall i,
\label{eq:cap}
\end{equation}
where $b_{\min}, b_{\max}$ are the minimum and maximum battery capacity within which the battery should be operating. The power consumed by the battery at time $i$ is denoted as
\begin{equation}
f(x_i)= \frac{[x_i]^+}{h \eta_{\text{ch}}} - \frac{\eta_{\text{dis}}[x_i]^-}{h}=\frac{\max(0,x_i)}{h\eta_{\text{ch}}} - \frac{\eta_{\text{dis}}\max(0,-x_i)}{h}, 
\label{finverse}
\end{equation}
where $x_i$ must lie in the range from $X_{\min}=\delta_{\min}h$ to $X_{\max}={\delta_{\max}h}$.
The ramping constraint is given as
\begin{equation}
    x_i \in [X_{\min}, X_{\max}].
    \label{eq:ramp}
\end{equation}

The battery is interfaced via an inverter. Let the inverter's efficiency be denoted as $\eta_{\text{inv}} \in (0,1]$. The modified battery charging and discharging efficiency is denoted as
\begin{gather}
    \eta_{\text{ch}}^* = \eta_{\text{ch}}\eta_{\text{inv}},~~~~
    \eta_{\text{dis}}^* = \eta_{\text{dis}}\eta_{\text{inv}}.
\end{gather}


\subsection{Inter-regional energy arbitrage model}
In this section, we formulate an optimal energy arbitrage model for energy storage participating in inter-regional energy arbitrage in two regions, shown as grid A and B in Fig. \ref{topo_fig}. 
Since both grids A and B have separate energy markets, the price levels for consumption or injection vary, this is because the size of the interconnector connecting the two grids is small compared to cumulative power needs in grid A and B\footnote{Interconnectors 
operated under a flow-based market would lead to a reduction in price fluctuations implying inter-regional energy trades would be marginally less profitable with their growth \cite{leuven2015cross}.}.
$P_i^{\text{\sc b,a}}$ and $P_i^{\text{\sc s,a}}$ denote the buying (consumption) and selling (injection) prices of electricity in grid A for time instant $i$, respectively. Similarly, $P_i^{\text{\sc b,b}}$ and $P_i^{\text{\sc s,b}}$ denote the buying and selling prices of electricity in grid B for time instant $i$, respectively.

As shown in Fig. \ref{topo_fig}, the battery is located in grid A. The efficiency of AC or DC interconnector linking grid A and B is given as $\eta_{\text{line}} \in (0,1]$\footnote{This implies $x$ MWh of energy supplied at one end of the interconnector will decrease to $\eta_{\text{line}} x$ MWh at the other end.}.
In this work, we assume that the energy storage and the interconnector are owned by different entities, thus, the interconnector owner levies a rent of $\zeta_i > 0, \forall i$ on the energy storage. The unit of $\zeta_i$ is in euros per MWh.

The buying price seen by the battery in grid A for charging from grid B is given as
\begin{equation}
    \widetilde{P}_i^{\text{\sc b,b}} = \frac{P_i^{\text{\sc b,b}} + \zeta_i}{\eta_{\text{line}}}.
    \label{eqbuypriceredef}
\end{equation}
The selling price seen by the battery in grid A for discharging from grid B is given as
\begin{equation}
    \widetilde{P}_i^{\text{\sc s,b}} = {(P_i^{\text{\sc s,b}} - \zeta_i)}{\eta_{\text{line}}}.
    \label{eqsellpriceredef}
\end{equation}
Note from \eqref{eqbuypriceredef} that the adjusted buying price seen in grid A increases due to interconnector line losses and interconnector rent, while the adjusted selling price denoted in \eqref{eqsellpriceredef} decreases due to interconnector line losses and interconnector rent.

\subsubsection{K1: Using either Grid A or B and not both}
This case uses either one of the power markets and not both simultaneously.
The objective function for the battery in this case is given as
\begin{equation}
   \min \sum \Big\{ \min (P_i^{\text{\sc b,a}}, {\widetilde{P}_{i}^{\text{\sc b,b}}})\frac{[x_i]^+}{\eta_{\text{ch}}} - \max(P_i^{\text{\sc s,a}}, \widetilde{P}_{i}^{\text{\sc s,b}})[x_i]^-\eta_{\text{dis}} \Big\}
   \label{eqobjcase1}
\end{equation}
\eqref{eqobjcase1} can be solved using linear programming formulation proposed in \cite{hashmi2019optimal, hashmi2019optimization}.
This case assumes that both grids A and B are fully available to charge or discharge, which may not be true.


\subsubsection{K2: Multi-market Charging and Discharging}
Cost functions for operating the battery using grids A and B are given by
\begin{align}
    & C_i^{\text{\sc a}} = P_i^{\text{\sc b,a}} \frac{[x_{i}^{\text{ \sc a}}]^+}{\eta_{\text{ch}}^*} - P_i^{\text{\sc s,a}} [x_{i}^{\text{ \sc a}}]^- \eta^*_{\text{dis}}, \\ 
    & C_i^{\text{\sc b}} = \widetilde{P}_{i}^{\text{\sc b,b}} \frac{[x_{i}^{\text{ \sc b}}]^+}{\eta^*_{\text{ch}}} - \widetilde{P}_{i}^{\text{\sc s,b}} [x_{i}^{\text{ \sc b}}]^- \eta^*_{\text{dis}}, \ \text{respectively},
\end{align}
where $x_{i}^{\text{ \sc a}}, x_{i}^{\text{ \sc b}} \in [X_{\min}, X_{\max}]$ denote the change in battery charge level due to grid A and B, respectively.
Thus, the battery ramping constraint is redefined as
\begin{equation}
    (x_{i}^{\text{ \sc a}} + x_{i}^{\text{ \sc b}}) \in [X_{\min}, X_{\max}].
    \label{eq:ramp2a}
\end{equation}
 {Note that it is crucial that $x_{i}^{\text{ \sc a}}$ and $x_{i}^{\text{ \sc b}}$ have the same sign, i.e., the battery can either charge or discharge at any given moment, but not both simultaneously. This requirement can be modelled via the constraint,}
 \begin{equation}
     x_{i}^{\text{ \sc a}} \cdot x_{i}^{\text{ \sc b}} \geq 0 \quad \forall i.
     \label{eqbilinear}
 \end{equation}
 As it is, constraint \eqref{eqbilinear} is intractable, as it involves a bilinear term and its feasible region is also non-convex. Conventional approaches for handling bilinear terms, such as the McCormick relaxation \cite{mccormick1976computability}, are not suitable as the relaxation can be very weak and provide impractical simultaneous charging and discharging solutions for the battery. Therefore, we instead choose to exactly reformulate the nonlinear constraint \eqref{eqbilinear} into a set of linear constraints, via the following crucial observation. For every time period, $i$, in the space of variables $(x_{i}^{\text{ \sc a}}, x_{i}^{\text{ \sc b}})$, the feasible region is either the non-negative orthant ($x_{i}^{\text{ \sc a}} \geqslant 0, x_{i}^{\text{ \sc b}} \geqslant 0$) or the non-positive orthant ($x_{i}^{\text{ \sc a}} \leqslant 0, x_{i}^{\text{ \sc b}} \leqslant 0$).  Since this feasible region is piecewise convex \cite{nagarajan2019adaptive, yang2022optimal}, we can equivalently reformulate it to capture the disjunctive union of these two regions by introducing two binary variables $z_i^{\text{ch}}, z_i^{\text{dis}} \in \{0,1\}$, and by satisfying the following inequalities $\forall i$: 
 \begin{subequations}
 \begin{align}
     & x_{i}^{\text{ \sc a}} \geqslant z_i^{\text{ch}} \cdot X_{\min}, \quad x_{i}^{\text{ \sc a}} \leqslant z_i^{\text{dis}} \cdot X_{\max},\\
     & x_{i}^{\text{ \sc b}} \geqslant z_i^{ch} \cdot X_{\min}, \quad x_{i}^{\text{ \sc b}} \leqslant z_i^{\text{dis}} \cdot X_{\max},\\
     & z_i^{\text{ch}} + z_i^{\text{dis}} = 1.0. 
 \end{align}
 \label{eq:linearization}
 \end{subequations}
 
 \vspace{-15pt}

\subsection{Epigraph-based mathematical formulation}

{The inter-regional energy arbitrage problem discussed so far {for K2} is in the form of a constrained minimization problem of a convex piecewise linear cost function with linear constraints and binary variables. This is given as:}

\begin{topbot}
    \text{Original problem ($P_{\text{orig}}$)}
\end{topbot}
\vspace{-10pt}
\begin{IEEEeqnarray}{ l C l }
    \label{eq:original_formulation}
    \text{Objective function:~~} 
     &\min{ \sum_i^N (C_i^{\text{\sc a}} ~ + ~C_i^{\text{\sc b}}) },& \label{eqoriOBJ} \\
    \text{subject to:} & \eqref{eq:cap}, \eqref{eq:ramp2a}, \eqref{eqbilinear} & \nonumber \\
    \hline \nonumber
\end{IEEEeqnarray}

\pagebreak

{A problem of this form can be transformed into an equivalent MILP formulation by forming the epigraph problem as in:}

\begin{topbot}
    \text{Epigraph formulation ($P_{\text{MILP}}$)}
\end{topbot}
{\allowdisplaybreaks
\begin{IEEEeqnarray}{ l C l }
    \label{eq:epi_formulation}
    \text{Objective function:} &~& \nonumber \\
    \min{ \quad \Big (~\sum_i t_{i}^{\text{ \sc a}} + \sum_i t_{i}^{\text{ \sc b}} ~\Big)}, &~& \label{eqoriOBJ_epi} \\
    \text{subject to:} && \nonumber \\
    \text{Cost function segment 1 and 2:} &~& \nonumber \\
    P_i^{\text{\sc b,a}} \frac{x_{i}^{\text{ \sc a}}}{\eta^*_{\text{ch}}} \leq t_{i}^{\text{ \sc a}}, &~& \forall~~ i,
    \label{eq14}\\
    P_i^{\text{\sc s,a}} x_{i}^{\text{ \sc a}} \eta^*_{\text{dis}} \leq t_{i}^{\text{ \sc a}}, &~& \forall~~ i, \label{eq15}\\
    \text{Cost function segment 3 and 4:} &~& \nonumber \\
    \widetilde{P}_{i}^{\text{\sc b,b}} \frac{x_{i}^{\text{ \sc b}}}{\eta^*_{\text{ch}}} \leq t_{i}^{\text{ \sc b}}, &~& \forall~~ i, \label{eq16}\\
    \widetilde{P}_{i}^{\text{\sc s,b}} x_{i}^{\text{ \sc b}} \eta^*_{\text{dis}} \leq t_{i}^{\text{ \sc b}}, &~& \forall~~ i, \label{eq17}\\   
 \text{Ramp constraint:~} &~& \nonumber \\
 x_{i}^{\text{ \sc a}} \in [X_{\min}, X_{\max}],&~& ~\forall~ i,\\
 x_{i}^{\text{ \sc b}} \in [X_{\min}, X_{\max}],&~& ~\forall~ i, 
 \label{eq:rampxB} \\
 (x_{i}^{\text{ \sc a}} + x_{i}^{\text{ \sc b}}) \in [X_{\min}, X_{\max}],&~& ~\forall~ i, 
 \label{eq:ramp2} \\
 \text{Capacity constraint:~} &~& \nonumber \\
 \sum_{j=1}^i \{x_{j}^{\text{ \sc a}} + x_{j}^{\text{ \sc b}}\} \leq b_{\max}-b_0,~&~& \forall~ i, \label{eq20} \label{eqcap1}\\
 -\sum_{j=1}^i \{x_{j}^{\text{ \sc a}} + x_{j}^{\text{ \sc b}}\} \leq b_0- b_{\min},~&~& \forall~ i, \label{eq21} \label{eqcap2} \\
 \text{Constraints} \ \eqref{eq:linearization} \ \text{to ensure:~} x_{i}^{\text{ \sc a}}\cdot x_{i}^{\text{ \sc b}}\geq 0, ~\forall i &~&   \\
    \hline \nonumber
\end{IEEEeqnarray}
}

{Re-formulating the problem in this manner is beneficial, as $P_{\text{MILP}}$ is composed of two epigraphs for the two regions, for each time instance. This decouples the storage participating in two energy markets simultaneously. Thus, this formulation can be easily modified for more than two energy markets simultaneously, provided \eqref{eqbilinear} is extended for the desired number of market zones.}
However, the focus of this work would be to assess two energy markets. 
In the next section, we extend $P_{\text{MILP}}$ to consider interconnector flows and offshore wind injections while considering the interconnector capacity.

\pagebreak

\section{Interconnector flow \& Capacity Blocking }
\label{section4}
AC or DC interconnections linking two separate energy markets have a limited capacity, denoted as $L_{\max}$. 
But the inter-regional energy arbitrage formulation, $P_{\text{MILP}}$, assumes that the AC or DC interconnector is not limiting the battery charging and discharging, as shown in \eqref{eq:rampxB}. Thus, $P_{\text{MILP}}$ assumes $\max(|X_{\min}|,X_{\max}) \leq L_{\max}$ and interconnector capacity is assumed to be available for all time instances.
This may not be true, as the interconnectors are often used for many other purposes (renewable energy curtailment reduction, security of supply, energy trade, ancillary services \cite{kaushal2019overview} etc.) and may not be fully available. 
In this section, we develop \textit{operating envelopes} for battery optimization based on the interconnector flows, and the offshore wind injections (OWI).

\subsection{Impact of interconnector flow: test case 1}
We denote the flow in the NEMO link as $L^{\text{\sc ab}}_i$, where $L^{\text{\sc ab}}_i>0$ implies power injected from Belgium to the UK and vice versa.
The power flow in the interconnector might affect the battery's ability to charge or discharge.
Due to the flow, 
\eqref{eq:rampxB} needs to be updated considering the interconnector flow. The operating envelopes are defined as follows:

\begin{equation}
        X_{\max}^{\text{adj}} = \begin{cases}
                       \max\big(0, \min(X_{\max}, L_{\max} + L^{\text{\sc ab}}_i)\big), \text{~if $L^{\text{\sc ab}}_i < 0$}, \\
                        X_{\max}, \text{~if $L^{\text{\sc ab}}_i \geq 0$}.
                    \end{cases}
         \label{eqlimmax}           
\end{equation}
\vspace{-10pt}
\begin{equation}
        X_{\min}^{\text{adj}} = \begin{cases}
                       X_{\min}, \text{~if $L^{\text{\sc ab}}_i < 0$}, \\
                       \min\big(0, \max(X_{\min}, -L_{\max} + L^{\text{\sc ab}}_i)\big), \text{~if $L^{\text{\sc ab}}_i \geq 0$}.
                    \end{cases}
         \label{eqlimmax}           
\end{equation}

\eqref{eq:rampxB} thus gets modified into
\begin{equation}
     x_{i}^{\text{ \sc b}} \in [X_{\min}^{\text{adj}}, X_{\max}^{\text{adj}}], ~\forall~ i.
\end{equation}


\subsection{Impact of interconnector: test case 2}
\begin{figure}[!htbp]
	\center
	\includegraphics[width=6.6in]{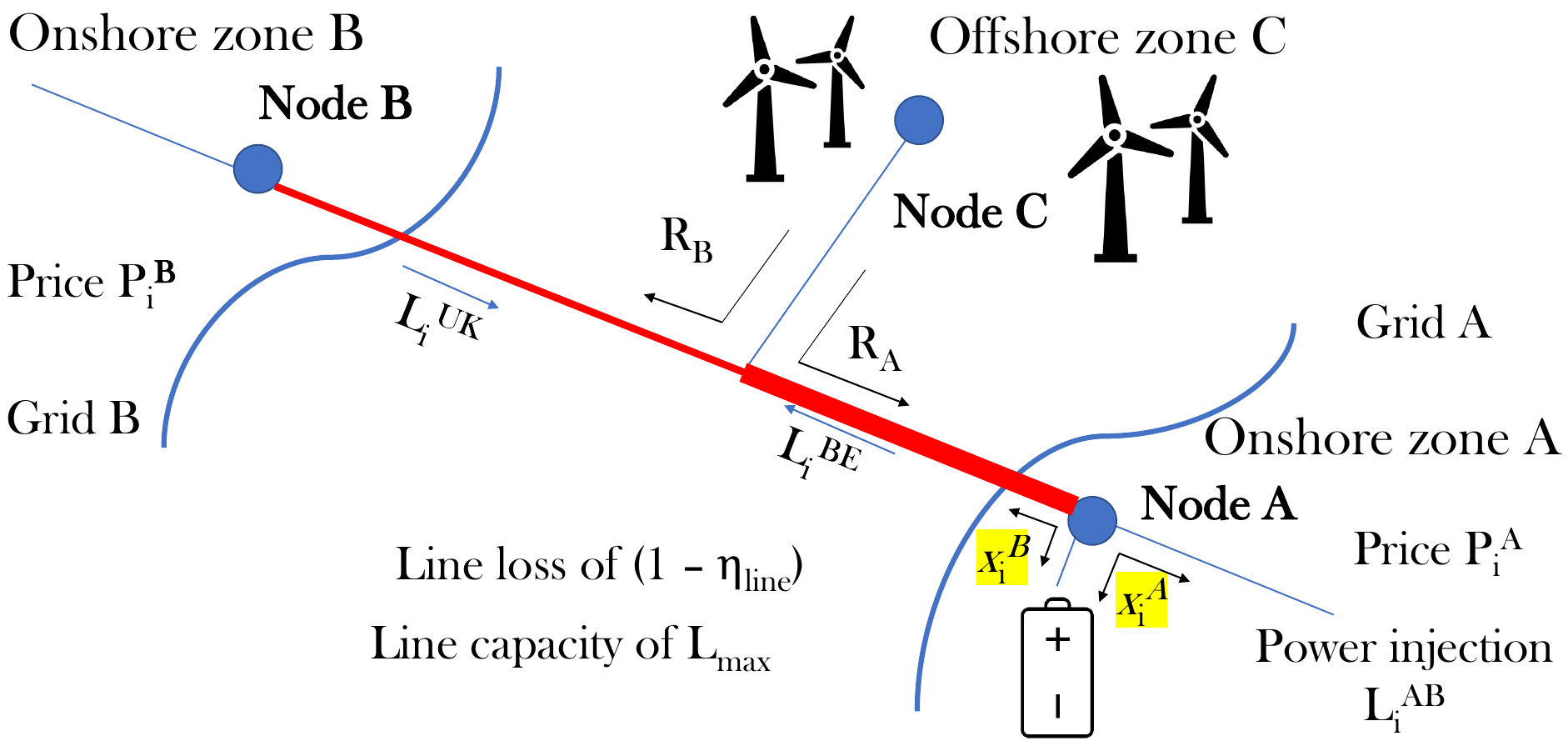}
	\vspace{-1pt}
	\caption{Interconnector (shown in red) with offshore injection. The associated variables are indicated. 
 }
	\label{fig:casegen}
\end{figure}

The Belgian and the UK transfer capacities to the PEEI have a capacity of 3.5 and 1.4 GW respectively. The different line limits are denoted as $L_{\max}^{\text{\sc be}}$ and $L_{\max}^{\text{\sc uk}}$. The OWPP at the PEEI affects the flows in the lines towards Belgium and the UK. 
The flow in these lines are denoted as $L^{\text{\sc be}}_i$ and $L^{\text{\sc uk}}_i$.
Fig. \ref{fig:casegen} shows the stylized diagram showcasing the interconnector flow and OWPP.

The operating envelope of the interconnector including the OWPP is the intersection of the operating regions for Belgium and the UK.
The operating region at time $i$ in the Belgian portion of the interconnector is given as
\begin{equation}
        X_{\max}^{\text{adj, BE}} = \begin{cases}
                       \max\big(0, \min(X_{\max}, L_{\max}^{\text{\sc be}} + L^{\text{\sc be}}_i)\big), \text{~if $L^{\text{\sc be}}_i < 0$}, \\
                        X_{\max}, \text{~if $L^{\text{\sc be}}_i \geq 0$}.
                    \end{cases}
         \label{eqlimmax}           
\end{equation}
\vspace{-10pt}
\begin{equation}
        X_{\min}^{\text{adj, BE}} = \begin{cases}
                       X_{\min}, \text{~if $L^{\text{\sc be}}_i < 0$}, \\
                       \min\big(0, \max(X_{\min}, -L_{\max}^{\text{\sc be}} + L^{\text{\sc be}}_i)\big), \text{~if $L^{\text{\sc be}}_i \geq 0$}.
                    \end{cases}
         \label{eqlimmax}           
\end{equation}

The operating region at time $i$ in the UK portion of the interconnector is given as
\begin{equation}
        X_{\max}^{\text{adj, UK}} = \begin{cases}
                       \max\big(0, \min(X_{\max}, L_{\max}^{\text{\sc uk}} + L^{\text{\sc uk}}_i)\big), \text{~if $L^{\text{\sc uk}}_i < 0$}, \\
                        X_{\max}, \text{~if $L^{\text{\sc uk}}_i \geq 0$}.
                    \end{cases}
         \label{eqlimmax}           
\end{equation}
\vspace{-10pt}
\begin{equation}
        X_{\min}^{\text{adj, UK}} = \begin{cases}
                       X_{\min}, \text{~if $L^{\text{\sc uk}}_i < 0$}, \\
                       \min\big(0, \max(X_{\min}, -L_{\max}^{\text{\sc uk}} + L^{\text{\sc uk}}_i)\big), \text{~if $L^{\text{\sc uk}}_i \geq 0$}.
                    \end{cases}
         \label{eqlimmax}           
\end{equation}

\eqref{eq:rampxB} in  $P_{\text{MILP}}$ thus gets modified into
\begin{equation}
     x_{i}^{\text{ \sc b}} \in [X_{\min}^{\text{adj, \sc be}}, X_{\max}^{\text{adj, \sc be}}]  \cap
     [X_{\min}^{\text{adj, \sc uk}}, X_{\max}^{\text{adj, \sc uk}}]
     , ~\forall~ i.
\end{equation}

\subsection{Capacity blocking for emergency services}

\begin{figure}[!htbp]
	\center
	\includegraphics[width=2.99in]{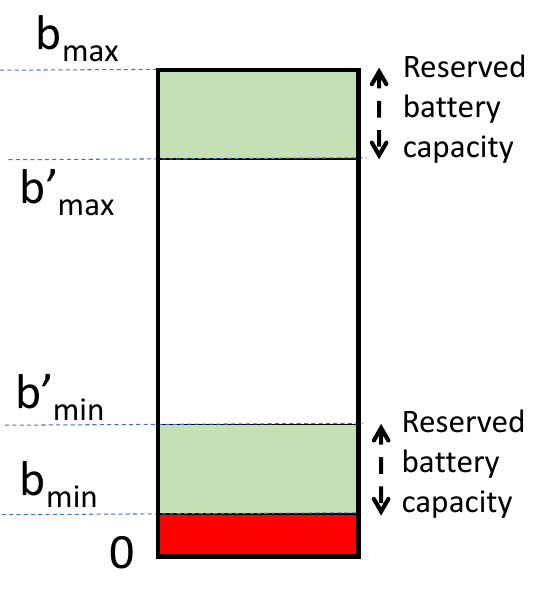}
	\vspace{-1pt}
	\caption{Capacity blocking for emergency services.}
	\label{fig:stack}
\end{figure}

Based on the use cases detailed in Section \ref{sec:usecases}, we model the capacity blocking as shown in Fig. \ref{fig:stack}.
Capacity blocking typically involves reserving a certain portion of the energy storage system's capacity for specific use cases, such as emergency services. This means that a predetermined amount of the battery's energy would be kept in reserve to ensure that critical services have access to power during emergencies or grid disruptions.
  
For  $P_{\text{MILP}}$, capacity blocking for emergency services are modeled by replacing 
$b_{\max}$ 
and $b_{\min}$ 
by $b_{\max}^{'}$ and $b_{\min}^{'}$, 
respectively, in 
\eqref{eqcap1} and \eqref{eqcap2}.

The total blocked battery is denoted as 
\begin{equation}
    b_{\text{block}} = (b_{\max} - b_{\max}^{'}) + (b_{\min}^{'} - b_{\min}).
\end{equation}

\pagebreak

\section{Performance indices}
\label{section5}
The performance indices, detailed in the following subsections, are used to evaluate the numerical simulations.

In this work, we define the interconnector utilization factor in percentage is defined as 
\begin{equation}
    \texttt{UF}= 100 \times \frac{\sum_{i=1}^T |L^{\text{\sc ab}}_i|}{L_{\text{lim}} T}, i \in \{1,..,T\}.
\end{equation}

\subsection{Revenue from arbitrage}
The revenue from performing energy arbitrage in multiple markets is denoted as
$
    R = -\sum_i^N \{C_i^{\text{\sc a}} ~ + ~C_i^{\text{\sc b}}\}.
$

\subsection{Cycles of operation of the battery}
The life of the battery is often defined in terms of cycle life and calendar life. The cycle life is determined by the operational cycles of charging and discharging. The relationship between the cycle of operation and depth of discharge is not a linear one. In this work, we utilize Algorithm 1 in \cite{hashmi2018long} for calculating the cycles of operation.
The input to this algorithm is the charge and discharge trajectory for the whole year calculated by solving $P_{\text{MILP}}$.
The cycle counting of a battery is crucial for assessing the financial viability of installing a battery.

\subsection{Simple payback period}
A simple payback period is the number of years required to recover the initial investment. The interest and operational cost is not considered. The revenue for 1 year of simulation is used for calculating a simple payback period denoted as
\begin{equation}
    \text{SPP} = \frac{\text{Investment cost}}{\text{1 year of revenue}}.
\end{equation}

\subsection{Tuning optimal value blocking level}
Finding the ideal battery capacity that can be blocked is a critical question. We present two models for calculating the ideal value of $b_{\text{block}}$.

\subsubsection{M1: Knee point}
The knee point of the plot of SPP versus $b_{\text{block}}$ {(Fig. \ref{fig:stack1})} is {calculated using \cite{kneepoint}}.

\subsubsection{M2: Based on the calendar life}
{For Fig. \ref{fig:stack1}}, the capacity corresponding to the calendar life of the battery is selected as the desired value of $b_{\text{block}}$.

\subsection{Computation time}
For the simulations, the computation time is calculated using 1000 Monte Carlo simulation runs. 
Simulations are performed on HP Intel(R) Core(TM) i7 CPU, 1.90GHz, 32 GB RAM personal computer on Matlab 2021a.

\section{Numerical case study}
\label{section6}
The two case studies first introduced in Section \ref{section2} are now presented. The battery characteristics considered in these case studies are detailed in Tab. \ref{tab:batparameters}.

\begin{table}[!htbp]
\centering
\caption{Battery parameter used for simulations}
\begin{tabular}{l|l} 
\hline
Attributes & Value \\
\hline \hline
Cost                         & 100 euros/kWh  \\ 
Rated capacity ($b_{\max}$)              & 1 MWh          \\ 
Minimum operational capacity ($b_{\min}$)  & 0.1 MWh        \\ 
Max charging rate    ($\delta_{\max}$)         & 0.5 MW         \\ 
Min discharging rate   ($\delta_{\min}$)      & -0.5 MW        \\ 
Charging efficiency      ($\eta_{\text{ch}}$)    & 0.95           \\ 
Discharging efficiency ($\eta_{\text{dis}}$)      & 0.95           \\ 
Converter efficiency    ($\eta_{\text{conv}}$)     & 0.95           \\ 
Initial charge level  ($b_0$)       & 0.5 MWh        \\ 
Cycle life (100\% DoD)       & 7200           \\ 
Calendar life                & 10 years       \\
\hline
\end{tabular}
\label{tab:batparameters}
\end{table}

As mentioned, a data-driven market model is used based on the years 2019 and 2020. It is observed, however, that there are missing dates in the price and interconnector flow time series. {To clean the data, any day with 2 or more consecutive hours of missing data (i.e. 2 samples) is ignored from the analysis.} Single missing hours are replaced by the average of the previous and subsequent values. Applying this rule, we received 357 days of data for 2019 and 331 days for 2020.
The suboptimal model presented in \eqref{eqobjcase1} is not evaluated, as the optimality gap observed between K1 denoted in \eqref{eqobjcase1} generates 24.2\% less revenue compared to $P_{\text{MILP}}$ model. 
Further, we cannot consider interconnector flows in \eqref{eqobjcase1}.


\begin{figure}[!htbp]
     \centering
    \begin{subfigure}{0.89\textwidth}
        \raisebox{-\height}{\includegraphics[width=\textwidth]{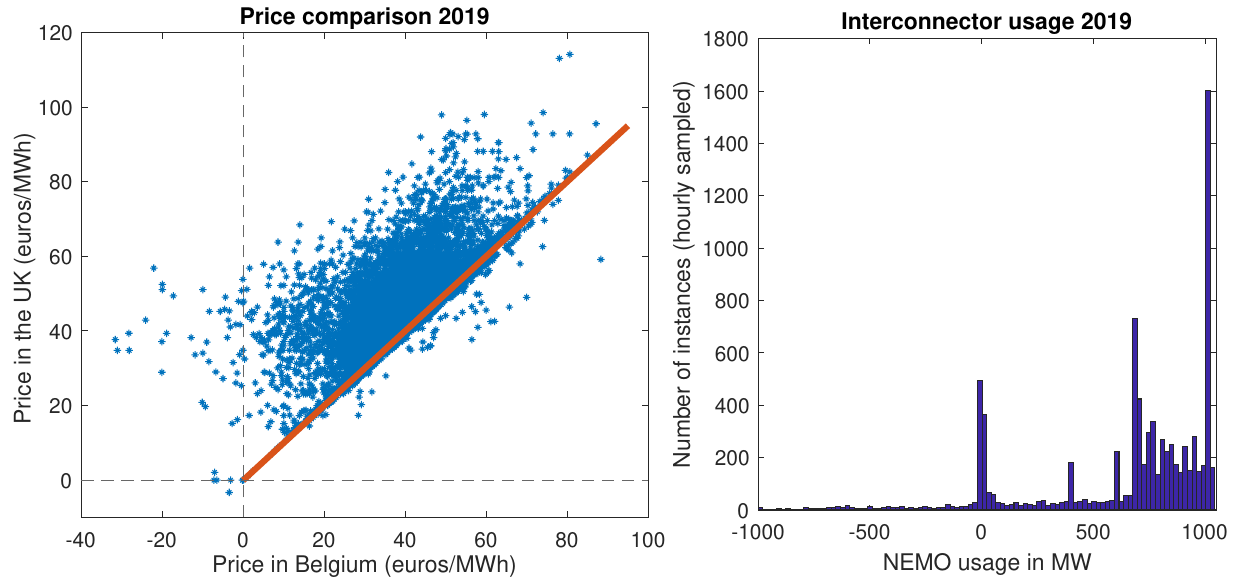}}
        \caption{Price comparison and NEMO link utilization for 2019}
    \end{subfigure}
    \hfill
    \begin{subfigure}{0.89\textwidth}
        \raisebox{-\height}{\includegraphics[width=\textwidth]{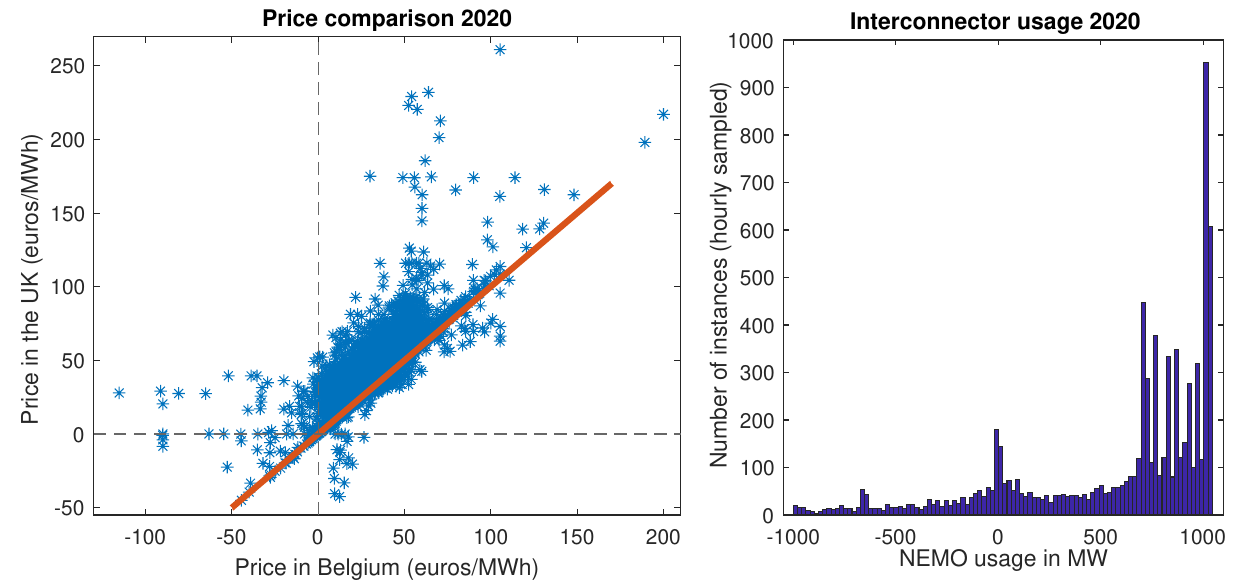}}
        \caption{Price comparison and NEMO link utilization for 2020}
    \end{subfigure}
    \caption{\small{Comparing the electricity prices for Belgium and the UK for the years 2019 and 2020. The mean utilization factor for the years 2019 and 2020 are 68.56 and 68.34\%. The convention used: positive value means injection from Belgium to the UK}}
    \label{fig:price}
\end{figure}

Note that the arbitrage mechanism proposed in this work relies on the convexity of the problem. The formulated optimization is no longer convex for negative electricity price levels. However, negative prices are not frequent. For 2019, the electricity prices were negative for 50 and 1 hourly instances for Belgium and the UK respectively. This consists of 0.3\% of the total horizon.
For 2020, 136 and 91 hourly instances of negative prices were observed for Belgium and the UK respectively, consisting of 1.4\%.
Since the negative prices are not frequent, the price level used for evaluating the numerical simulation is saturated at the lower limit of zero. The mean electricity prices for 2019 are 37.9 and 47.24 euros/MWh for Belgium and the UK\footnote{The electricity prices in the UK are in pounds/MWh. This is converted to euros/MWh. The conversion factor used in the work is fixed at 1.16.}.
The mean electricity prices for 2020 are 31.28 and 40.79 euros per MWh for Belgium and the UK.
{For traceability of numerical results, we assume that $\zeta_i$ is constant for all time instances.}

Fig. \ref{fig:price} compares the electricity prices in Belgium and the UK for the years 2019 and 2020. Note that pre- and post-BREXIT, the electricity prices have similar trends, where the prices in the UK are more often higher than in Belgium. This is evident from the fact that most price instances are higher than the red line. The red line denotes the hypothetical line where prices in the two regions are equal. This skewness in the electricity price motivates us to explore the possibility of utilizing energy storage batteries in the two energy markets.
The observation made in Fig. \ref{fig:price} 
is that 
the interconnector is predominantly used to inject power from Belgium and extract it in the UK. 
The electricity price distribution for the years 2019 and 2020 is similar.

{
The MILP-based inter-regional energy arbitrage,  $P_{\text{MILP}}$, code
described in this paper, is publicly available at
\href{https://github.com/umar-hashmi/Inter-Refional-Energy-Arbitrage}{{github.com/umar-hashmi/Inter-Refional-Energy-Arbitrage}} \cite{gitcode}.
}


\pagebreak

\subsection{Case study 1: NEMO Link}
The NEMO link connects mainland Belgium with the UK grid via the North Sea. The undersea cable length is 140 km long. 
The capacity of the line is 1 GW, with line losses in the order of 2.5\% \cite{linknemoloss}. 
The battery is assumed to be placed at the Belgian end of the NEMO link. {We have two objectives in performing this case study:}
\begin{enumerate}
    \item A sensitivity analysis on the interconnector rent towards performance indices, and
    \item Selecting $b_{\text{block}}$ for emergency services: Evaluating M1, M2. 
\end{enumerate}

{To accomplish this, we consider three scenarios:}
\begin{itemize}
   \item { \textbf{C1 (Lower bound on benefits)}: considers storage participating in only the local Belgian grid.}
    \item { \textbf{C2 (Upper bound on benefits)}: considers the storage participating in both the Belgian and the UK energy markets without any constraint on charging and discharging due to interconnector unavailability.}
    \item {  \textbf{C3}: considers the storage participating in both grids, while also considering the interconnector capacity constraints. The battery revenue for C3 should be lower than or at best equal to C2.}
\end{itemize}

\subsubsection{\textbf{Objective 1: Sensitivity towards interconnector rent}}
\begin{figure}[!htbp]
	\center
	\includegraphics[width=6.9in]{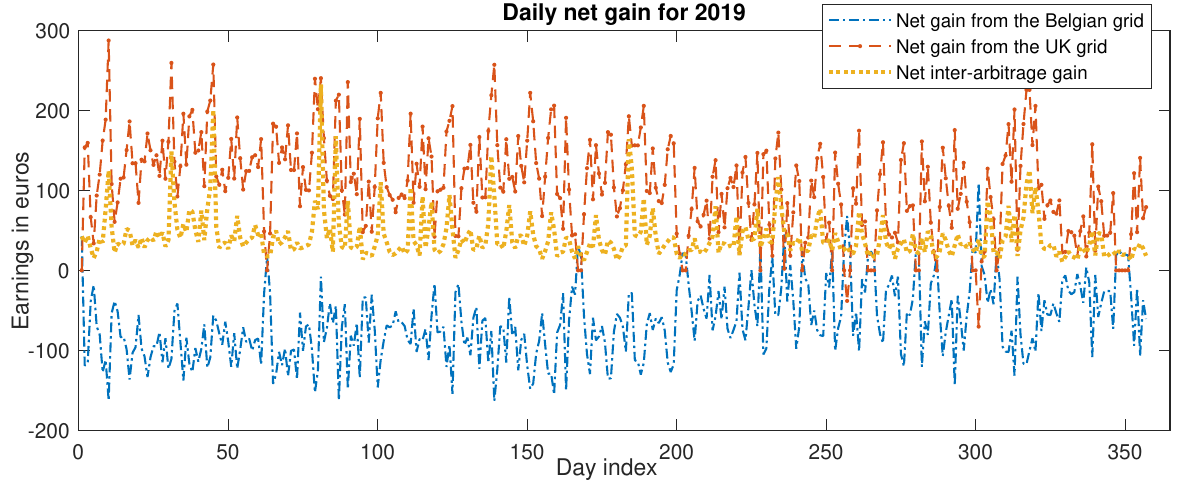}
	\vspace{-1pt}
	\caption{Revenue separation for inter-regional arbitrage for the year 2019 without interconnector flows.}
	\label{fig:div}
\end{figure}

The interconnector auctioned the capacity and levied a charge for its utilization proportional to the power level \cite{linkAUCTION}. 
Fig. \ref{fig:div} shows the separated revenue for inter-regional arbitrage. Note that the battery utilizes the Belgian grid for charging and almost entirely discharges into the UK grid. This is expected due to the price difference bias between these regions also highlighted in Fig. \ref{fig:price}.

\begin{figure}[!htbp]
	\center
	\includegraphics[width=5.8in]{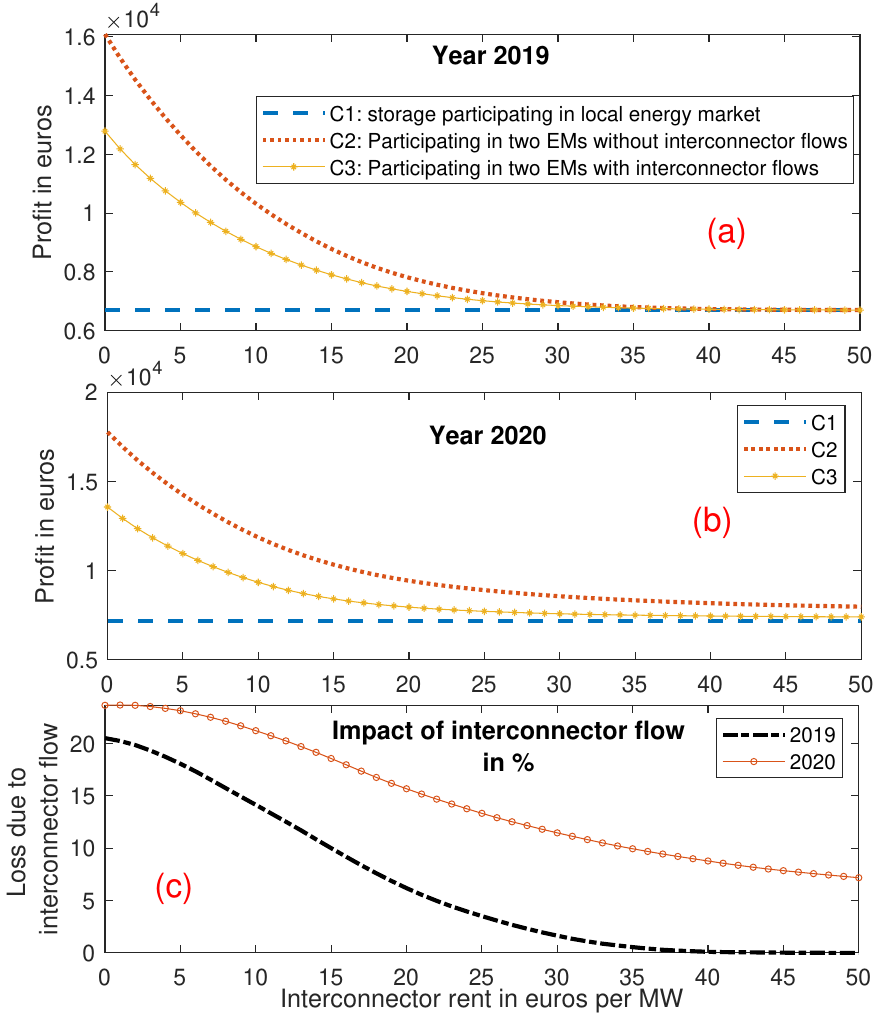}
	\vspace{-1pt}
	\caption{Profit due to battery with interconnector rent.}
	\label{fig:res1}
\end{figure}

\begin{figure}[!htbp]
	\center
	\includegraphics[width=5.8in]{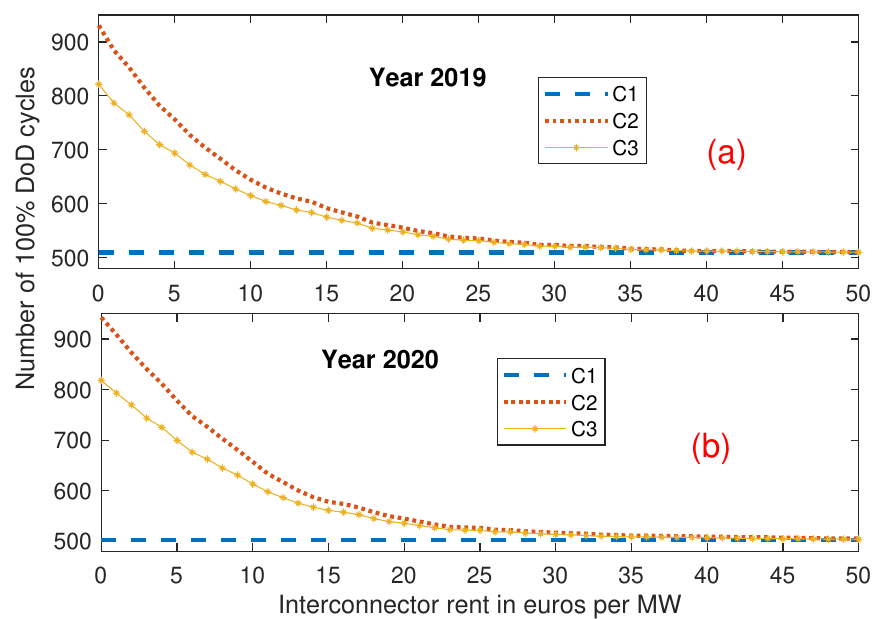}
	\vspace{-1pt}
	\caption{Cycles of operation with interconnector rent.}
	\label{fig:res2}
\end{figure}

\begin{figure}[!htbp]
	\center
	\includegraphics[width=5.7in]{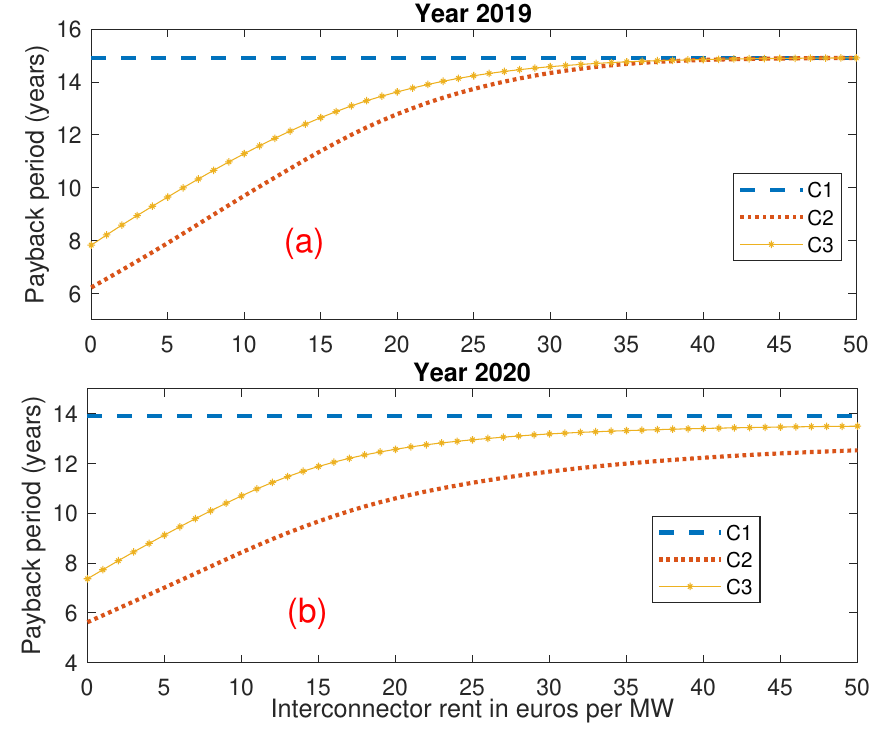}
	\vspace{-1pt}
	\caption{Simple payback period with interconnector rent.}
	\label{fig:res3}
\end{figure}

\begin{figure}[!htbp]
	\center
	\includegraphics[width=5.7in]{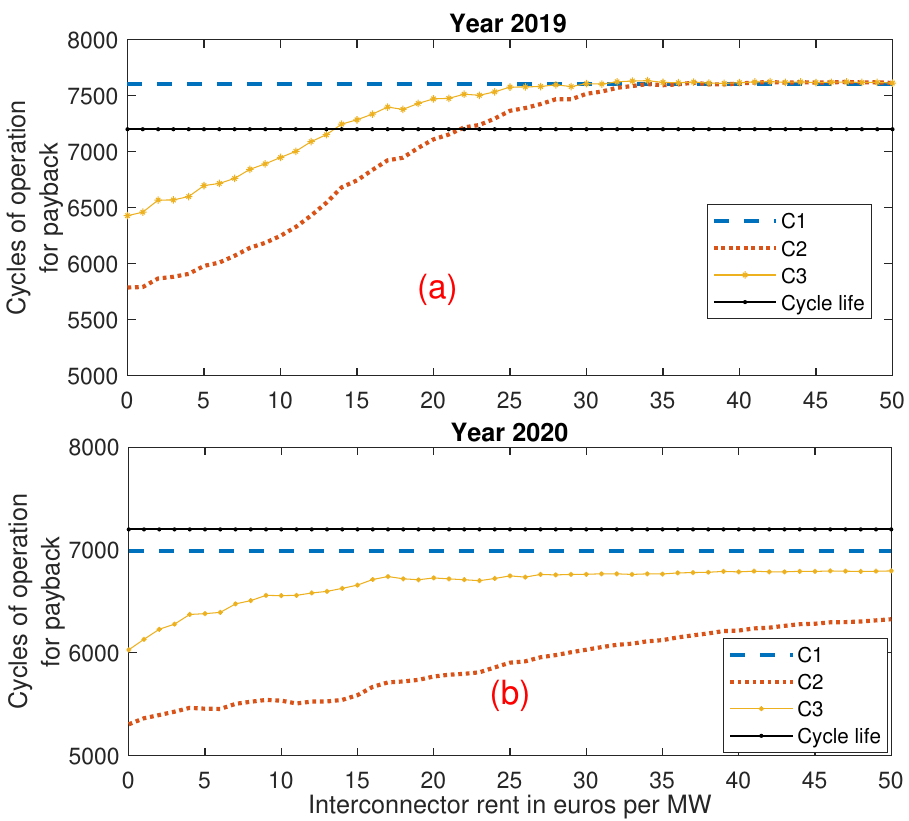}
	\vspace{-1pt}
	\caption{Cycles of operation for SPP with interconnector rent.}
	\label{fig:res4}
\end{figure}

Fig. \ref{fig:res1}, \ref{fig:res2}, \ref{fig:res3}, and \ref{fig:res4} shows the performance indices for models C1, C2 and C3 with different levels of interconnector rent depicted on the \textit{x}-axis.


Fig. \ref{fig:res1} (a) and (b) show the profit due to the battery for the scenarios C1, C2 and C3 for the years 2019 and 2020, respectively.
Note that model C1 is independent of interconnector rent, as in this case, the battery considers only the Belgian grid.
Also, C2 which does not consider interconnector flows generates a higher profit
compared to C3. However, the disparity in scenarios C2 and C3 decreases as the interconnector rent increases, thus, making it not profitable to use the interconnector at all.
Fig. \ref{fig:res1} (c) shows the impact of interconnector flows on the profitability of the battery. For small values of interconnector rent, the impact of {interconnector congestion} on profitability exceeds 20\% and 24\% for the years 2019 and 2020, respectively.

Fig. \ref{fig:res2} shows the total number of 100\% DoD cycles the battery performed for the years 2019 and 2020.
The cycles of operation for models C2 and C3 converge to C1 for high levels of interconnector rent, as it is no longer profitable to use the interconnector for energy exchange.

{Fig. \ref{fig:res3} illustrates the simple payback period computed based on profits from the years 2019 and 2020. Notably, for lower interconnector rent values, the simple payback period for the C3 case can be as short as 8 years, indicating the potential profitability of the battery system (as described in Table \ref{tab:batparameters}). It is important to emphasize that these calculations do not account for additional revenue streams that the battery could generate by participating in additional market products within both regions.
These findings hold promise for prospective battery owners. However, it is worth noting that as more batteries are integrated into the grid, the marginal value of installing a battery is expected to decline. This reduction in value is attributed to the decrease in electricity price volatility when there is an increased level of energy storage in the system, a phenomenon discussed in  \cite{hashmi2018effect}.}

Fig. \ref{fig:res4} plots the number of cycles the battery would require to reach its payback investment. The black line shows the cycle life of the battery.
Observe that for the year 2019, the battery inter-regional arbitrage is not profitable for model C3 for an interconnector rent of more than 13 euros/MW. For model C1, the battery is not profitable for the year 2019.
However, for the year 2020, model C1 is profitable.
This is due to the price volatility, which is higher in the year 2020 compared to 2019. The price variance for the year 2020 was 287.6 and 401.9 in Belgium and the UK. However, for the year 2019, it was 137.7 and 139.1 for Belgium and the UK; substantially less than in 2020.

\subsubsection{\textbf{Objective 2: Capacity blocking}}
In this part of the case study related to the NEMO link, we evaluate two models for selecting capacity blocking level, $b_{\text{block}}$. 
Tab. \ref{tab:capblock} lists the key outcomes for models M1 and M2 regarding the capacity blocked, SPP, and cycles of operation.

\begin{table}[!htbp]
\centering
\caption{Evaluating models for capacity blocking }
\begin{tblr}{
  cell{2}{1} = {r=4}{},
  cell{2}{2} = {r=2}{},
  cell{4}{2} = {r=2}{},
  cell{6}{1} = {r=4}{},
  cell{6}{2} = {r=2}{},
  cell{8}{2} = {r=2}{},
  vlines,
  hline{1-2,6,10} = {-}{},
  hline{3,5,7,9} = {3-5}{},
  hline{4,8} = {2-5}{},
}
Years & Model & Metric                                   & C2    & C3    \\ \hline
2019  & M1 (SPP = 10)   & $b_{\text{block}}$ & 0.27  & 0.04  \\
      &  (SPP = 10)     & 100\% DoD cycles         & 6677  & 6666  \\
      & M2 ($b_{\text{block}}$=0.4)   & Payback (years)                          & 11.26 & 13.64 \\
      &       & 100\% DoD cycles         & 6295  & 6780  \\
2020  & M1 (SPP = 10)   & $b_{\text{block}}$ & 0.41  & 0.12  \\
      &  (SPP = 10)     & 100\% DoD cycles         & 5993  & 6399  \\
      & M2  ($b_{\text{block}}$=0.4)  & Payback (years)                          & 9.81  & 12.86 \\
      &        & 100\% DoD cycles         & 6005  & 6676  
\end{tblr}
\label{tab:capblock}
\end{table}

\begin{figure}[!htbp]
	\center
	\includegraphics[width=5.59in]{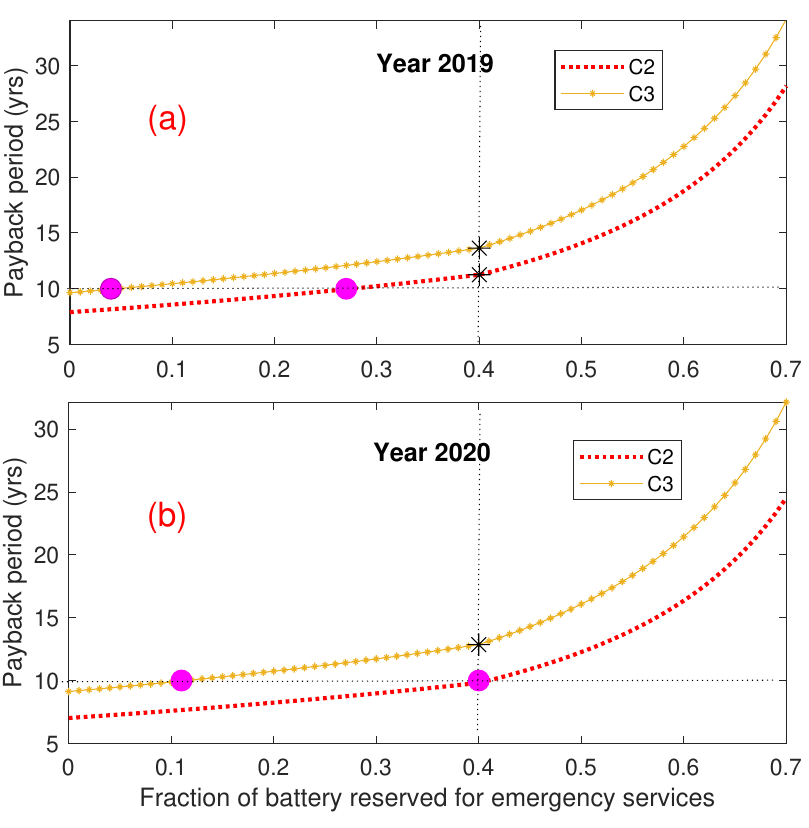}
	\vspace{-1pt}
	\caption{Profitability of battery considering capacity blocking. Pink circles denote M1 and black asterisk denote M2.}
	\label{fig:stack1}
\end{figure}

Fig. \ref{fig:stack1} plots the SPP and $b_{\text{block}}$ for C2 and C3 models. Note that due to limited battery life, the $b_{\text{block}}$ is lower for M1 compared to M2.
For the year 2019 and 2020, the $b_{\text{block}}$ for C3 model are 4 and 12\%.

\begin{figure}[!htbp]
	\center
	\includegraphics[width=5.29in]{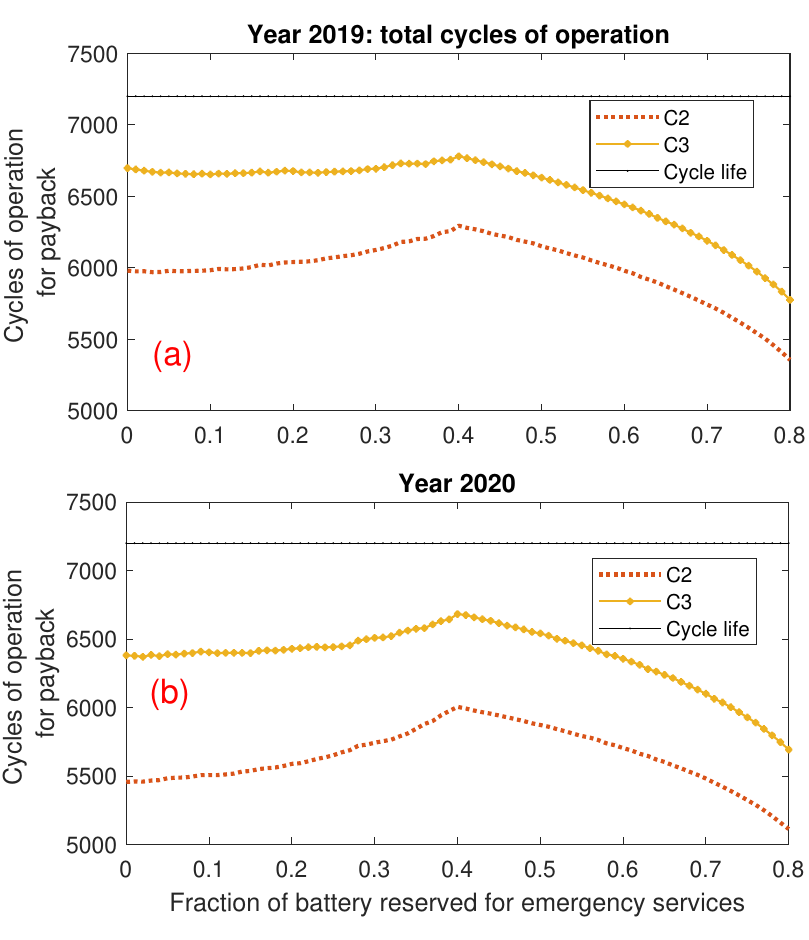}
	\vspace{-1pt}
	\caption{Cycles of operation due to inter-regional arbitrage comparison with capacity blocking}
	\label{fig:stack2}
\end{figure}

Fig. \ref{fig:stack2} plots the 100\% DoD cycles the battery performed versus $b_{\text{block}}$. 
Note that as the battery capacity is reserved for other emergency applications, the number of cycles the battery is performing due to inter-regional arbitrage reduces.

The key takeaways of this case study are:

$\bullet$ \textit{Interconnector Flow:} When evaluating the financial feasibility of energy storage engaged in inter-regional energy arbitrage, it is essential to account for interconnector flows. These flows can significantly affect potential gains.\\
$\bullet$ \textit{Interconnector Rent:} A rising level of interconnector rent can diminish the attractiveness of multi-region energy arbitrage as an investment opportunity.\\
$\bullet$ \textit{Profitability of Belgian-based Batteries:} Batteries positioned at the Belgian end of the NEMO link have the potential for profitability. However, 
this is contingent on 
factors such as price volatility and interconnector flow dynamics.\\
$\bullet$ \textit{Battery Lifecycle Consideration:} Evaluating the financial viability of battery installation should incorporate a thorough assessment of both the battery cycle and calendar life.\\
$\bullet$ \textit{Capacity Blocking for Emergency Services:} When implementing capacity blocking strategies for emergency services, it is imperative to factor in the battery's cycle and calendar life to ensure sustained profitability.





\pagebreak

\subsection{Case study 2}
{The second case study examines interconnectors connecting the PEEI in the North Sea. It is assumed that inter-regional markets operate under a flow-based market design as described in 
Sec. \ref{section2}. Under such a case, the opportunities for inter-regional arbitrage are significantly reduced as interconnector flow becomes saturated in a single direction.}

\begin{figure}[!htbp]
	\center
	\includegraphics[width=7in]{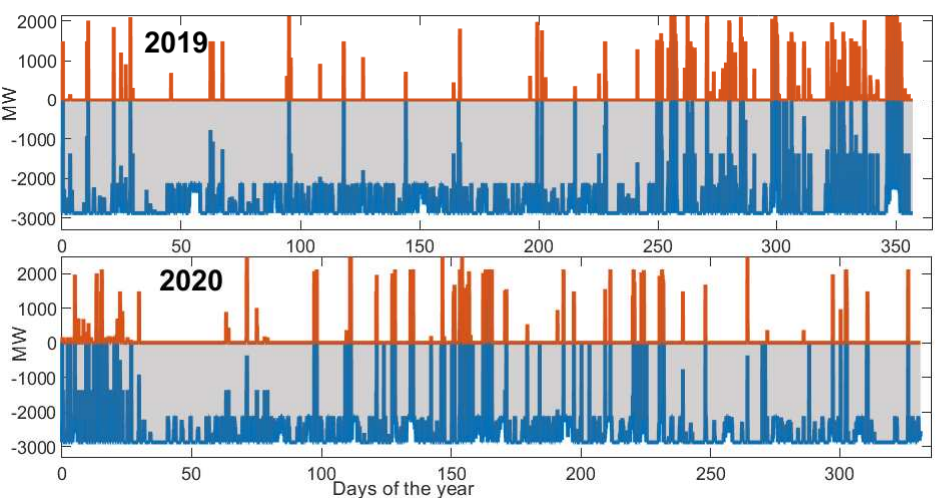}
	\vspace{-3pt}
	\caption{Operating envelopes, i.e., feasible range of interconnector capacity available for inter-regional arbitrage.}
	\label{fig:feasiblerange}
\end{figure}

{As illustrated in Fig. \ref{fig:feasiblerange}, it is generally possible for the battery to charge from the UK grid but discharging into the UK grid for the majority of the time is not possible without previously blocked capacity. Given this constraint under the flow-based market design, battery owners must participate in an interconnector capacity auction to reserve interconnector capacity apriori. As such, in this case study, we assume the interconnector rent is awarded at a fixed rate of 5 euros per megawatt-hour (MWh).}

\begin{figure}[!htbp]
	\center
	\includegraphics[width=5.5in]{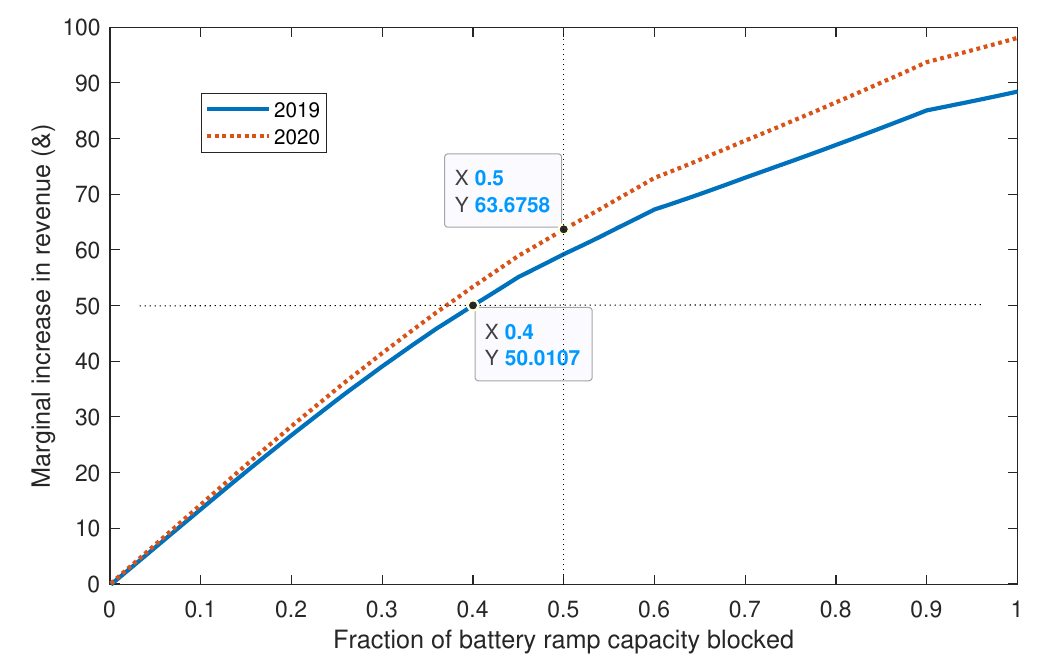}
	\vspace{-1pt}
	\caption{Marginal increase in the revenue with the fraction of interconnector capacity reserved compared to battery ramp rate.}
	\label{fig:case22}
\end{figure}

{Fig. \ref{fig:case22} presents the marginal increase in storage revenue when compared to participating exclusively in the Belgian energy market. It is evident from the graph that the marginal revenue increase is nearly linear up to 50\% of the ramping power limit. However, it diminishes beyond that point. This analysis underscores that storage owners can significantly boost their revenue, surpassing an 88\% increase, by reserving interconnector capacity equal to the battery's ramping limit. However, for a comprehensive optimal bidding strategy, further numerical evaluations are necessary to inform battery owners' decisions.}

\begin{figure}[!htbp]
	\center
	\includegraphics[width=5.5in]{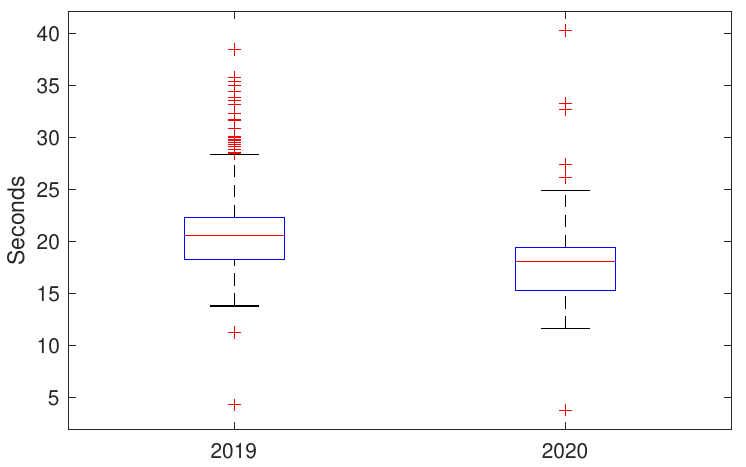}
	\vspace{-1pt}
	\caption{Computation time for 1000 MC for 2019, 2020}
	\label{fig:computationTime}
\end{figure}

Fig. \ref{fig:computationTime} shows the computational time distributions for 1000 Monte Carlo (MC) simulations for different levels of interconnector capacity blocked for the years 2019 and 2020.

The key takeaways of the case study are:\\
$\bullet$ When interconnector flows are calculated based on the flow-based market design, the interconnector saturates in one direction, leaving no room for inter-regional energy arbitrage. Under this condition, the battery owners should bid to reserve interconnector capacity. For interconnector capacity reserved, we observe an increase in battery revenue by more than 88\% compared to the battery participating only in the local grid.\\
$\bullet$ The median computation time of 1000 MC simulations is less than 21 seconds for one whole year.

\pagebreak

\section{Conclusions}
\label{section7}

This paper explores the evolving landscape of the future power system, characterized by a growing web of AC and DC interconnectors that bridge previously separate energy markets. The interconnectors serve as a vital link between these distinct energy markets. The authors propose an advanced model for optimal multi-market energy storage arbitrage, augmenting it with emergency services that can be extended to system operators.

The proposed modelling includes
battery ramping and capacity constraints, and the utilization of operating envelopes derived from interconnector capacity, efficiency, dynamic energy injection, and offshore wind generation forecasts for the day ahead. 
The key findings of this paper are as follows: Firstly, it {develops} a Mixed-Integer Linear Programming (MILP) model for energy storage, enabling simultaneous energy arbitrage across multiple markets while adhering to constraints like interconnector capacity, efficiency, and charge-discharge sequencing. 
The energy arbitrage model is designed to accommodate two distinct electricity price dynamics, with specialized considerations for bilinear terms.
Secondly, the paper introduces a profitability assessment framework for energy storage, accounting for a portion reserved for emergency services and the battery's limited lifespan. Lastly, the study conducts realistic case studies, including NEMO Link and interconnections between Belgium and the UK connected to the PEEI, exploring factors such as interconnector rent's impact on performance, the allocation of battery capacity for emergencies, and the potential reservation of interconnector capacity for inter-regional energy trading.

The paper introduces two different models for determining the upper limit of capacity that energy storage owners can allocate for emergency services. To validate their approach, the authors conduct numerical case studies,  focusing on the interconnections between Belgium and the UK.
In summary, we present a comprehensive model for optimizing grid-scale batteries in the context of interconnected energy markets, with a focus on both profitability and emergency services, demonstrated through real-world case studies.

\pagebreak

\section*{Acknowledgement}
This project has received funding from the CORDOBA project funded by Flanders Innovation and Entrepreneurship (VLAIO) in the framework of the spearhead cluster for blue growth in Flanders (Blue Cluster) – Grant number HBC.2020.2722. 

\pagebreak

\bibliographystyle{IEEEtran}
\bibliography{reference}

\onecolumn
\appendix

\subsection{Implementing MILP formulation for inter-regional energy arbitrage}
\label{appendix:lpmatrix}

The matrix format for the optimization problem $P_{\text{MILP}}$ is denoted as 
minimize
${f}^T X$, subject to ${A}X\leq b$, and $X \in [lb, ub]$.
The dimension of $A$ is 13Nx6N, $b$ is 13Nx1, $X$ and $f$ are of size 6Nx1, and N denotes the number of samples in the horizon of optimization.



\begin{equation}
f\text{=}{\begin{bmatrix}
	0\\
	:\\
	0\\
	0\\
	:\\
	0\\
	1\\
	:\\
	1\\
 1\\
	:\\
	1\\
    0\\
	:\\
	0\\
 0\\
	:\\
	0\\
	\end{bmatrix}},~~
X \text{=} {\begin{bmatrix}
	x_{1}^{\text{\sc a}}\\
	:\\
	x_{N}^{\text{\sc a}}\\
    x_{1}^{\text{\sc b}}\\
	:\\
	x_{N}^{\text{\sc b}}\\
	t_{1}^{\text{\sc a}}\\
	:\\
	t_{N}^{\text{\sc a}}\\
 	t_{1}^{\text{\sc b}}\\
	:\\
	t_{N}^{\text{\sc b}}\\
    z_{1}^{\text{\sc a}}\\
    :\\
    z_{N}^{\text{\sc a}}\\
    z_{1}^{\text{\sc b}}\\
    :\\
    z_{N}^{\text{\sc b}}\\
	\end{bmatrix}},~~
lb\text{=}
{\begin{bmatrix}
	X_{\min}\\
	:\\
	X_{\min}\\
 \hline
 	X_{\min}\\
	:\\
	X_{\min}\\
\hline
	T_{\min}^{\text{\sc a}}\\
	:\\
	T_{\min}^{\text{\sc a}}\\
 \hline
 	T_{\min}^{\text{\sc b}}\\
	:\\
	T_{\min}^{\text{\sc b}}\\
 \hline
	0\\
	:\\
	0\\
  \hline
	0\\
	:\\
	0\\
	\end{bmatrix}}
\leq
{\begin{bmatrix}
	x_{1}^{\text{\sc a}}\\
	:\\
	x_{N}^{\text{\sc a}}\\
 \hline
    x_{1}^{\text{\sc b}}\\
	:\\
	x_{N}^{\text{\sc b}}\\
 \hline
	t_{1}^{\text{\sc a}}\\
	:\\
	t_{N}^{\text{\sc a}}\\
 \hline
 	t_{1}^{\text{\sc b}}\\
	:\\
	t_{N}^{\text{\sc b}}\\
 \hline
    z_{1}^{\text{\sc a}}\\
    :\\
    z_{N}^{\text{\sc a}}\\
 \hline
    z_{1}^{\text{\sc b}}\\
    :\\
    z_{N}^{\text{\sc b}}\\
	\end{bmatrix}} \leq
ub\text{=}
{\begin{bmatrix}
	X_{\max}\\
	:\\
	X_{\max}\\
 \hline
	X_{\max}\\
	:\\
	X_{\max}\\
 \hline
	T_{\max}^{\text{\sc a}}\\
	:\\
	T_{\max}^{\text{\sc a}}\\
 \hline
	T_{\max}^{\text{\sc b}}\\
	:\\
	T_{\max}^{\text{\sc b}}\\
 \hline
1\\
    :\\
1\\
 \hline
1\\
    :\\
1\\
	\end{bmatrix}}, ~
	 {b} = {\begin{bmatrix}
	\text{zeros}(N,1) \\
 \text{zeros}(N,1) \\
 \text{zeros}(N,1) \\
 \text{zeros}(N,1) \\
 \hline
	b_{\max} - b_0\\
	:\\
	b_{\max} - b_0\\
	b_0- b_{\min}\\
	:\\
	b_0- b_{\min}\\
 \hline
  \text{ones}(N,1) X_{\max} \\
 -\text{ones}(N,1) X_{\min}\\
 \text{zeros}(4N,1) \\
 \text{ones}(N,1) \\
	\end{bmatrix}}.
\label{mateqsame}
\end{equation}
where $T_{\min}$ and $T_{\max}$ are bounds on $t_i$. Since these bounds are not known to us, we choose $T_{\min}$ to be negative with a large magnitude and $T_{\max}$ to be positive with a large magnitude.

{ \small{
\begin{equation}
{A} = { 
\renewcommand\arraystretch{1.3}
\mleft[ \begin{array}{c |c |c | c | c|c }
\texttt{diag}_N(P^{\text{\sc b,a}}_i/\eta^*_{\text{ch}}) & \texttt{zeros}_N & -1 \times \texttt{diag}_N & \texttt{zeros}_N & \texttt{zeros}_N & \texttt{zeros}_N \\ 
\texttt{diag}_N(P^{\text{\sc s,a}}_i\eta^*_{\text{dis}}) & \texttt{zeros}_N & -1 \times \texttt{diag}_N & \texttt{zeros}_N & \texttt{zeros}_N & \texttt{zeros}_N\\ 
\texttt{zeros}_N &\texttt{diag}_N(\widetilde{P}_i^{\text{\sc b,b}}/\eta^*_{\text{ch}}) & \texttt{zeros}_N &  -1 \times \texttt{diag}_N & \texttt{zeros}_N  & \texttt{zeros}_N\\ 
\texttt{zeros}_N &\texttt{diag}_N(\widetilde{P}_i^{\text{\sc s,b}}\eta^*_{\text{dis}}) & \texttt{zeros}_N &  -1 \times \texttt{diag}_N & \texttt{zeros}_N  & \texttt{zeros}_N\\ 
\texttt{trigL}_N & \texttt{trigL}_N & \texttt{zeros}_N& \texttt{zeros}_N & \texttt{zeros}_N & \texttt{zeros}_N \\ 
-1 \times \texttt{trigL}_N & -1 \times \texttt{trigL}_N & \texttt{zeros}_N& \texttt{zeros}_N & \texttt{zeros}_N & \texttt{zeros}_N\\ 
\texttt{diag}_N & \texttt{diag}_N & \texttt{zeros}_N& \texttt{zeros}_N & \texttt{zeros}_N & \texttt{zeros}_N \\ 
-1 \times \texttt{diag}_N & -1 \times \texttt{diag}_N & \texttt{zeros}_N& \texttt{zeros}_N & \texttt{zeros}_N & \texttt{zeros}_N\\ 
-1 \times \texttt{diag}_N & \texttt{zeros}_N & \texttt{zeros}_N & \texttt{zeros}_N & X_{\min} \times \texttt{diag}_N & \texttt{zeros}_N  \\
\texttt{diag}_N & \texttt{zeros}_N & \texttt{zeros}_N & \texttt{zeros}_N & \texttt{zeros}_N & - X_{\max} \times \texttt{diag}_N \\
\texttt{zeros}_N & -1 \times \texttt{diag}_N & \texttt{zeros}_N & \texttt{zeros}_N & X_{\min} \times \texttt{diag}_N & \texttt{zeros}_N  \\
\texttt{zeros}_N & \texttt{diag}_N  & \texttt{zeros}_N & \texttt{zeros}_N & \texttt{zeros}_N & - X_{\max} \times \texttt{diag}_N \\
\texttt{zeros}_N & \texttt{zeros}_N & \texttt{zeros}_N & \texttt{zeros}_N & \texttt{diag}_N & \texttt{diag}_N 
%
	\end{array} \mright]}
 \begin{array}{c }
\eqref{eq14} \vspace{3pt}\\
\eqref{eq15} \vspace{3pt}\\
\eqref{eq16} \vspace{3pt}\\
\eqref{eq17} \vspace{3pt}\\
\eqref{eq20} \vspace{3pt}\\
\eqref{eq21} \vspace{3pt}\\
\eqref{eq:ramp2} \vspace{3pt}\\
\eqref{eq:ramp2} \vspace{3pt}\\
\eqref{eq:linearization} \vspace{3pt}\\
\eqref{eq:linearization} \vspace{3pt}\\
\eqref{eq:linearization} \vspace{3pt}\\
\eqref{eq:linearization} \vspace{3pt}\\
\eqref{eq:linearization}
	\end{array}, 
\end{equation} }
}

where $\texttt{zeros}_N$ denotes a zero square matrix with $N$ rows and columns,
$\texttt{diag}_N$ denotes a diagonal matrix with 1 in the diagonal and zeros otherwise, 
$\texttt{trigL}_N$ denotes the lower triangular matrix of size $N$,
$\texttt{ones}_N$ denotes a square matrix of order $N$ with 1s and
$\texttt{diag}_N(J_i)$ denotes a diagonal matrix with $J_i$ in the diagonal and zeros otherwise.

Using a mixed integer linear programming solver, we solve the following problem $P_{\text{MILLP}}$.

\begin{topbot}
    \text{Implementing ($P_{\text{MILLP}}$)}
\end{topbot}
{\allowdisplaybreaks
\begin{IEEEeqnarray}{ l C  }
    \label{eq:epi_formulation}
    \text{Objective function:} & \min_X f^T X \nonumber \\
    \text{subject to:} & \text{(i) }A.X \leq b, \text{ (ii) }X \in [lb, ub], \text{ (iii) } X_{\text{binary}} \subset X \in \{0,1\}. \nonumber \\
    \hline \nonumber
\end{IEEEeqnarray}
}

In this work, we use MATLAB's \texttt{intlinprog} \cite{intlinprog} function for implementing $P_{\text{MILP}}$.

\end{document}